\documentclass[a4paper,12pt]{article}

\usepackage{amsmath, amssymb, amsthm}
\usepackage{mathtools}
\usepackage{enumerate}
\usepackage[dvips,final]{graphicx}
\usepackage[dvips,final]{epsfig}
\usepackage{caption}
\usepackage{subcaption}
\usepackage{enumitem}

\usepackage{hyperref}
\hypersetup{
	colorlinks=true,
	linkcolor=blue,
	urlcolor=cyan,
	citecolor=blue,
}

\usepackage{xcolor}

\definecolor{oran}{rgb}{.75,.4,0}
\definecolor{vio}{rgb}{0.4,0,0.6}
\definecolor{gre}{rgb}{0.1,0.6,0}
\definecolor{ora}{rgb}{0.9,0.7,0}
\definecolor{blgr}{rgb}{0,0.7,0.4}
\definecolor{brw}{rgb}{.5,0.25,0}

\definecolor{gray}{rgb}{0.8,0.8,0.8}
\definecolor{lgray}{rgb}{0.88,0.88,0.88}
\definecolor{llgray}{rgb}{0.95,0.95,0.95}

\usepackage{titlesec}
\titleformat{\section}{\bfseries}{\thesection}{1em}{}
\titleformat{\subsection}{\bfseries}{\thesubsection}{1em}{}
\titleformat{\subsubsection}{\itshape}{\thesubsubsection}{1em}{}

\numberwithin{equation}{section}

\newfont{\ctv}{msam10}

\newcommand{\bbox}{\mbox{\ctv \symbol{4}}}
\def\QED{{${}\hfill\bbox$}}
\newenvironment{pf}[1]{\par\vskip1mm{\noindent\it #1.}\ }{\QED\par
\vskip2mm}

\def\bpf{\begin{pf}}
\def\epf{\end{pf}}

\def\gpd{\hat g_{\delta}}
\def\bfe{\boldsymbol{e}}

\def\expe{\hbox{\rm e}}
\def\ve{\varepsilon}

\def\dd{\,\mathrm{d}}
\def\dive{\mathrm{\,div\,}}

\def\for{\mathrm{\ for\ }}
\def\ale{\mathrm{\ a.\,e.}}

\def\on{^{(n)}}
\def\sobo{W^{1,2}_0(\Omega)}

\def\vved{v^{\ve,\delta}}
\def\VV{\mathcal{V}_f}
\def\hgo{\hat\Gamma_0}
\def\hgom{\hat\Gamma_0^{-1}}

\def\som{\sin\,\omega}
\def\com{\cos\,\omega}

\def\scal#1{\left\langle #1 \right\rangle}

\def\real{\mathbb{R}}
\def\nat{\mathbb{N}}

\def\io{\int_{\Omega}}

\def\be{\begin{equation}\label}
\def\ee{\end{equation}}

\def\ber{\begin{eqnarray}}
\def\eer{\end{eqnarray}}

\def\bers{\begin{eqnarray*}}
\def\eers{\end{eqnarray*}}

\def\bpf{\begin{pf}}
\def\epf{\end{pf}}

\newtheorem{theorem}{Theorem}[section]
\newtheorem{lemma}[theorem]{Lemma}

\newtheorem{hypo}[theorem]{Hypothesis}
\newtheorem{propo}[theorem]{Proposition}

\begin{document}

\title{Analytic description of the moving moisture front in soils \thanks{The support from the European Union's Horizon Europe research and innovation programme under the Marie Sk\l odowska-Curie grant agreement No 101102708 and from the GA\v CR project 24-10586S is gratefully acknowledged.}
}

\author{
Bettina Detmann\thanks{University of Duisburg-Essen, Faculty of Engineering, Department of Civil Engineering, 45117 Essen, Germany, E-mail: {\tt bettina.detmann@uni-due.de}.}
\and Chiara Gavioli\thanks{Faculty of Civil Engineering, Czech Technical University, Th\'akurova 7, 16629 Praha 6 and Institute of Mathematics, Czech Academy of Sciences, \v Zitn\'a 25, 11567 Praha 1, Czech Republic, E-mail: {\tt gavioli@math.cas.cz}.}
\and Pavel Krej\v c\'{\i}\thanks{Faculty of Civil Engineering, Czech Technical University, Th\'akurova 7, 16629 Praha 6 and Institute of Mathematics, Czech Academy of Sciences, \v Zitn\'a 25, 11567 Praha 1, Czech Republic, E-mail: {\tt pavel.krejci@cvut.cz}.}
\and
Yanyan Zhang\thanks{School of Mathematical Sciences, Key Laboratory of MEA (Ministry of Education) and Shanghai Key Laboratory of PMMP, East China Normal University, Shanghai 200241, China, E-mail: {\tt yyzhang@math.ecnu.edu.cn}.}}

\date{}

\maketitle

\begin{abstract}
The fact that moisture propagates in soils at a finite speed is confirmed by natural everyday experience as well as by controlled laboratory tests. In this text, we rigorously derive analytical upper bounds for the speed of moisture front propagation under gravity for the solution to the Richards equation with compactly supported initial data. The main result is an explicit criterion describing a competition between gravity and capillarity, where the dominant effect is determined by the characteristics of the soil. If capillarity prevails, the initially wet regions remain wet for all times, while if gravity is dominant, moisture travels downward at a speed that is asymptotically bounded from below and above. As a by-product, we prove the existence and uniqueness of a solution to an initial value problem for the degenerate Richards equation on the whole space. Numerical simulations based on the proposed model confirm the theoretical predictions, with results that closely match experimental observations.

\bigskip


\medskip

\end{abstract}


\section{Introduction}\label{sec:intro}

\begin{figure}[htb]
	\begin{center}
		\includegraphics[height=5cm]{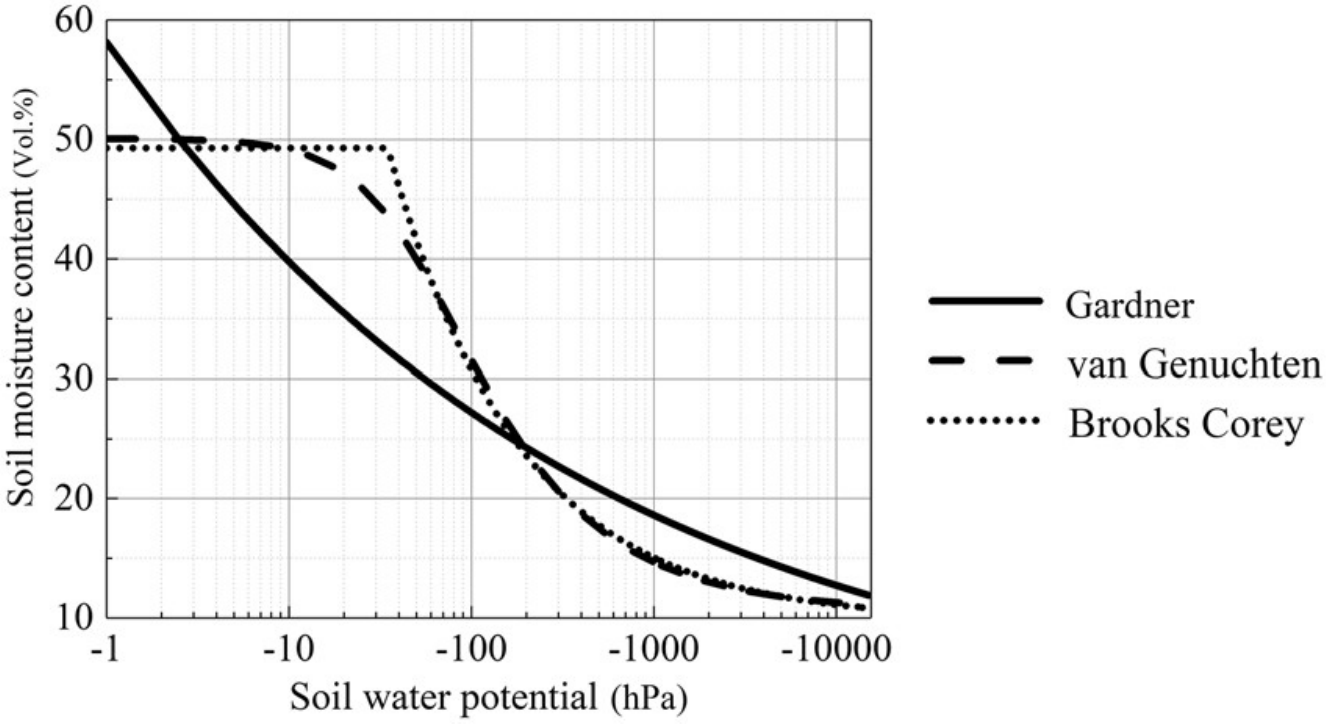}\qquad
		\includegraphics[height=5.1cm]{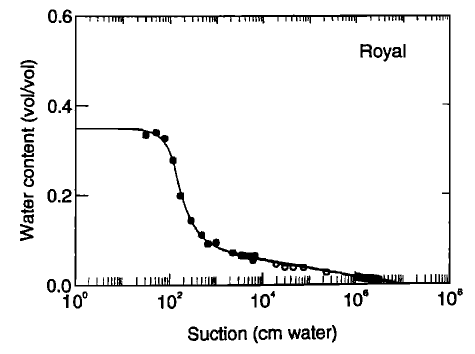}
		\caption{Comparison of water retention curves ($\theta(h)$). Left: Models with unbounded suction, where $\theta$ only asymptotically approaches zero as $h \to -\infty$. Right: Model with bounded suction (Rossi-Nimmo type), where $\theta$ reaches zero at a finite pressure head $h_*$.}
		\label{gf1}
	\end{center}
\end{figure}

The transport of water in soils is a fundamental problem in environmental physics and hydrology, with critical relevance to phenomena such as infiltration, groundwater recharge, and contaminant transport. This process is classically described by the Richards equation, a highly nonlinear, degenerate parabolic equation derived from Darcy's law and the conservation of mass:
\be{e0}
\theta_t - \dive \kappa(\theta) (\nabla h + \bfe_N) = 0,
\ee
with space variable $x \in \real^N$, $N \in \nat$, time variable $t \in (0,T)$, and unit vector $\bfe_N$ pointing against the gravity direction. The unknown functions are the pressure head $h = h(x,t)$ (negative in the unsaturated zone) and the volumetric water content $\theta = \theta(x,t) \in [0,n)$ (where $n$ denotes the porosity), while $\kappa$ is a given permeability function. We couple Eq.~\eqref{e0} with the Rossi-Nimmo soil water retention model (see \cite{rossi}), which provides a good description of water retention in both the dry and wet regions. This approach is an alternative to the asymptotic behavior of standard formulas like van Genuchten, whose reliance on infinite suction makes them unsuitable for applications involving low water contents, see Figure~\ref{gf1} (left) taken from \cite{pan}. The idea there is to introduce a finite cutoff pressure head $h_*$, where the liquid film continuity breaks and the retention curve vanishes, see Figure~\ref{gf1} (right) taken from \cite{rossi}. This value of $h_*$ corresponds to oven dryness and is of the order of $-10^5$ to $-10^7$ cm according to \cite{rossi}. In this framework, it is convenient to introduce a new unknown function: the shifted pressure head $u = u(x,t) \in [0,-h_*]$, which is defined as
\be{htou}
u = h-h_*.
\ee
The water retention function is assumed in the form
\be{e0a}
\theta = g(u),
\ee
where $g:\real\to\real$ is a continuous increasing function. In our shifted variable formulation, this implies the degeneracy condition $g(0) = 0$, see Hypothesis~\ref{hfg1}. Fluid diffusion in soil is then the result of competition between the fluxes $G \rho_w \nabla u$ driven by pressure and $-G \rho_w \bfe_N$ driven by gravity, where $\rho_w$ is the water mass density, and $G$ is the gravity acceleration constant.

\begin{figure}[htb]
	\centering
	\includegraphics[width=11cm]{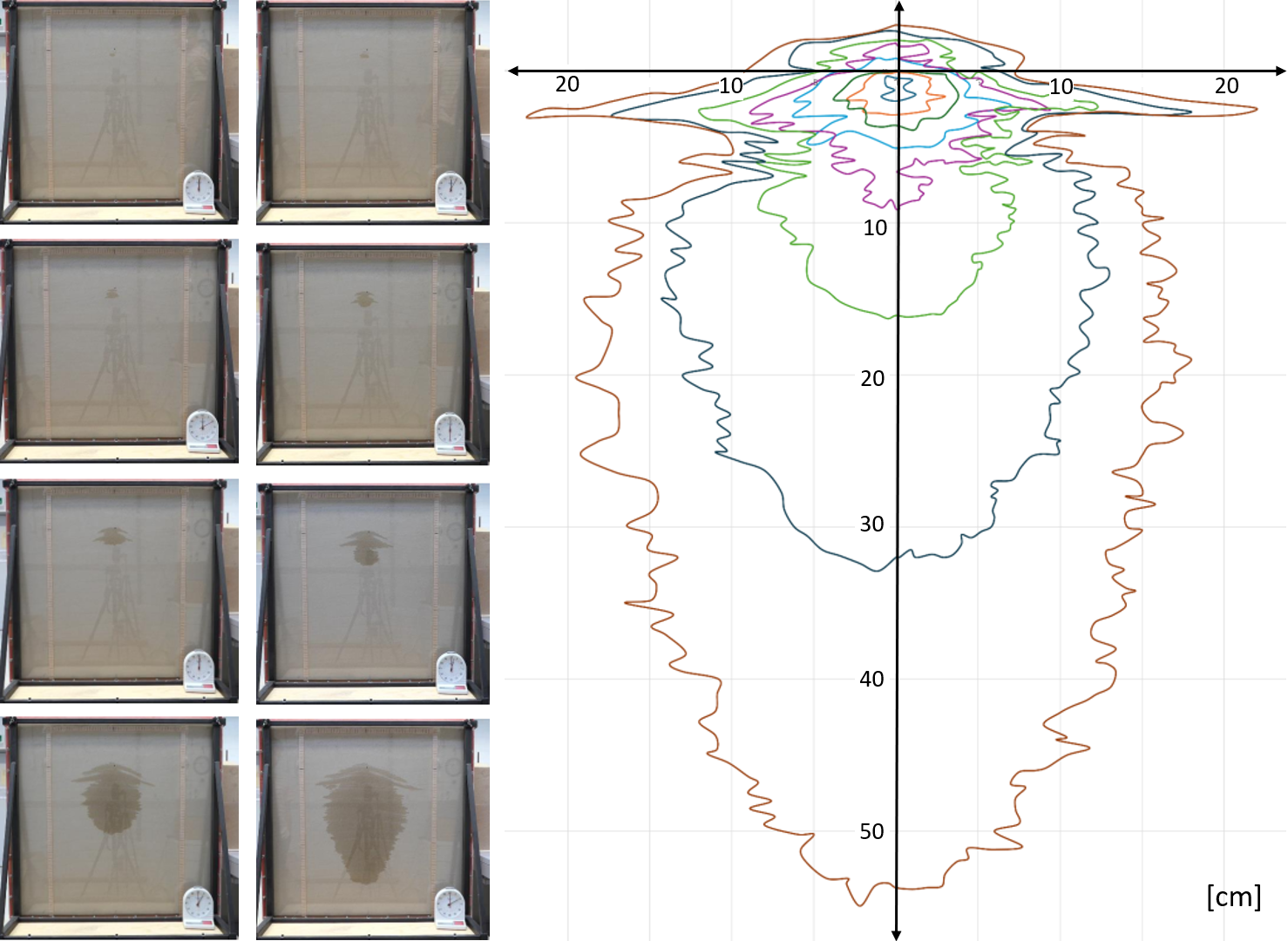}
	\caption{Experimental observation of two-dimensional wetting bulb propagation. Left: Time series photographs illustrating the development of the wetting bulb in initially dry sand, shown at eight specific times (1s, 5s, 10s, 30s, 60s, 120s, 300s, 600s). Right: Penetration figures of the wetting front corresponding to the times shown on the left. The water supply is in the origin of the coordinate system.}
	\label{fig:labres}
\end{figure}

\begin{figure}[htb]
\noindent
\begin{centering}
\includegraphics[width=9cm]{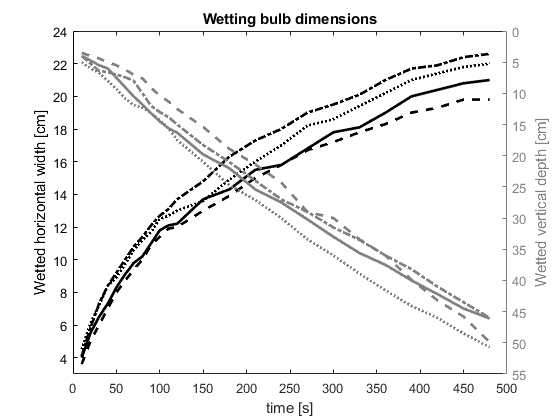}
\par\end{centering}
\caption{Evolution of the wetting bulb dimensions over time from four experimental tests. The $x$-axis is time in seconds. The wetted horizontal width (left axis, black) and the wetted vertical depth (right axis, gray) are shown in cm. Note the distinct propagation rates (linear vs.\ square-root growth) that motivate the subsequent analysis.}
\label{fig:bulbdim}
\end{figure}

The primary focus of this paper is the rigorous analytical investigation of the propagation properties and asymptotic behavior of solutions to Problem~\eqref{e0}--\eqref{e0a}. Specifically, we analyze the influence of the functions $\kappa$ and $g$ on the speed and shape of the wetting front. The theoretical investigation is motivated by laboratory experiments studying the two-dimensional propagation of a wetting bulb in initially dry sand, supplied by a point source. Such two-dimensional flow behavior is of critical importance in irrigation engineering for determining root water availability \cite{dong}; however, distinct experimental tests investigating two-dimensional flow remain relatively rare in the literature (see e.g., \cite{nat,dirk,wu}) compared to the widely studied one-dimensional problem (e.g., \cite{det,hav}). The laboratory test qualitatively shows the formation of a distinct, egg-shaped wetting bulb whose dimensions are governed by the interplay of capillary diffusion and downward gravity-driven convection (see Figure~\ref{fig:labres}). Initially, near the source, the flow is dominated by capillary diffusion, resulting in an isotropically radial front; however, the immediate influence of gravity quickly breaks this symmetry, leading to the characteristic egg-shaped profile. Furthermore, the experiments show that the upward capillary rise stops after a finite time, achieving a stationary height. Note that the large extent of the wetting front in the horizontal direction at the level of the outlet is an artifact that arises because the tube of the water supply disturbs the homogeneous filling of the sand into the cube. Crucially, the experimental results indicate distinct growth laws for the dimensions of the bulb: the vertical extent grows linearly with time, while the horizontal extent grows sublinearly with time (see Figure~\ref{fig:bulbdim}).

Although the laboratory experiments involve a localized, continuous water supply in the sand cube, in Problem~\eqref{e0}--\eqref{e0a} source terms are neglected because the qualitative long-time regimes driven by $\kappa$ and $g$ cannot be observed when a source term dictates the evolution. Thus, analyzing the source-free problem allows us to isolate and rigorously characterize the asymptotic properties driven solely by the nonlinear diffusion and convection. On the other hand, experiments without sources are difficult to evaluate because the moisture content quickly converges to zero. Furthermore, to maintain mathematical tractability, we neglect hysteresis effects in the constitutive relation \eqref{e0a}. While the Richards equation including hysteresis has been investigated in both engineering and mathematical contexts (e.g., \cite{kac_hys,ro_hys,sch_hys}), an analytical characterization of the support of the solutions in the hysteretic setting similar to \cite{gk2} remains an open problem. Our primary focus here is on the advancing front, which is accurately described by the primary wetting curve.

The Richards equation is extensively used in the literature to describe moisture dynamics, with a strong emphasis on numerical investigations of wetting bulb dimensions for irrigation applications (e.g., \cite{bai24,bai,mon,more,par}). These studies, along with works dedicated to developing numerical solvers \cite{bonan,zhang}, rely on computational results to assess the influence of hydraulic parameters. In contrast to these purely numerical treatments, our work offers a rigorous analytical perspective rooted in the equation's structure, which falls into the class of doubly-degenerate parabolic equations with a first-order convective term. The existence of compactly supported solutions is a defining feature of degenerate parabolic equations, originating with the porous medium equation \cite{vazquez}. The asymptotic behavior of this support has been the subject of extensive mathematical study. For pure degenerate parabolic problems, the support typically expands as $\sqrt{t}$ (see, e.g., \cite{andreucci,diaz,watanabe}), whereas the paper \cite{vespri} on doubly-nonlinear systems suggests that the propagation regime can also depend on the spatial dimension, particularly for small times. On the other hand, for first-order quasilinear equations dominated by convection, the front propagates linearly with $t$ (see \cite{diazver}). We propose here a contribution to this discussion by rigorously analyzing the anisotropic and combined effects of diffusion and gravitational transport. In contrast to the Barenblatt-type approach typical of parabolic problems, we employ traveling wave supersolutions to capture the anisotropic effects of gravity. Traveling wave solutions for the Richards equation have been studied in \cite{vD,ele_hys,wit} in the one dimensional setting. Although this approach generally precludes the inclusion of a source term, thus reinforcing our choice to analyze the source-free problem, it allows us to compute an explicit upper bound for the solution's support directly from the equation's structure. This enables us to characterize the long-time behavior without relying on numerical discretization. We show that for small times $t>0$, the moisture front propagates at most as $\sqrt{t}$ in all directions, independently of the space dimension and soil characteristics; this indicates that diffusion initially dominates the flow. However, for large times, the anisotropy induced by gravity dictates the regime: the diffusive $\sqrt{t}$ rate is preserved only in directions perpendicular to gravity, while in the direction of gravity the propagation rate shifts to at most $t$ (convection-dominated). In the direction opposite to gravity the support remains bounded by a constant surface, confirming the physical arrest of capillary rise. After some initial time interval, the upward front may also start moving downward.

The structure of the paper is the following. In Section~\ref{sec:auxi} we consider Eq.~\eqref{e0} on a bounded domain $\Omega \subset \real^N$ with zero Dirichlet boundary condition and state the main results of this paper. In Section~\ref{sec:exis} we first prove properties such as the positivity and uniqueness of solutions. Existence is proved under additional assumptions on $\kappa$ and $g$. The method consists in introducing two positive regularizing parameters $\ve$ and $\delta$ and proving that solutions to the approximating problems converge to the solution to the original problem. The results of Section~\ref{sec:exis} are summarized in Theorem~\ref{tfg1}. In Section~\ref{sec:mois} we show that, under suitable assumptions on $\kappa$ and $g$, solutions with compactly supported initial data have at each time $t>0$ a compact support $\Omega_t \subset \Omega$ which propagates with different speeds in the gravity direction and in directions perpendicular to gravity. We also describe the shape of the moving wet region and its asymptotic behavior. The results of Section~\ref{sec:mois} are summarized in Theorem~\ref{3t1}. Finally, in Section~\ref{sec:nume} we include numerical simulations to illustrate the theory developed in Section~\ref{sec:mois} for some relevant choices of $\kappa$ and $g$.


\section{Statement of the problem}\label{sec:auxi}

We fix a bounded Lipschitzian domain $\Omega \subset \real^N$. Defining the function
\be{deff}
f(u) = \int_0^u \kappa(g(v))\dd v
\ee
and recalling \eqref{htou}--\eqref{e0a}, we can rewrite equation \eqref{e0} in terms of the variable $u$ as
\be{e0u}
g(u)_t - \dive f'(u) (\nabla u + \bfe_N) = 0.
\ee
We consider equation \eqref{e0u} with a given initial condition
\be{e2}
u(x,0) = u_0(x) \quad \mbox{in} \ \Omega
\ee
and, for the moment, with homogeneous Dirichlet boundary condition
\be{e1}
u = 0 \quad \mbox{on} \ \partial\Omega.
\ee
More general boundary conditions will be discussed in Section~\ref{sec:mois}.
We define the space
\be{vpdef}
\VV = \{u \in L^\infty(\Omega): f(u) \in W^{1,2}_0(\Omega)\}.
\ee
The weak formulation of the problem corresponds to finding $u \in L^2(0,T;\VV)$ such that $g(u)_t \in L^2(0,T; W^{-1,2}(\Omega))$ and the identity
\be{1e1}
\io \big(g(u)_t \phi + f'(u) (\nabla u + \bfe_N)\cdot \nabla \phi\big)\dd x = 0
\ee
holds for all test functions $\phi\in W^{1,2}_0(\Omega)$.

\begin{hypo}\label{hfg1}
The functions $f,g:[0,\infty) \to \real$ are continuous together with derivatives $g'$, $f'$, and $f''$, and such that $f(0)=f'(0)=g(0)=0$, $f'(z) >0$, $g'(z)>0$ for $z>0$, and
\be{fg1}
\int_0^1 \frac{|f''(z)|^2}{g'(z)}\dd z < \infty.
\ee
\end{hypo}

The structural assumptions in Hypothesis~\ref{hfg1} are motivated by physical models of soil water retention that account for a finite pressure at oven dryness, most notably the formulation by Rossi and Nimmo \cite{rossi}. While standard models (such as van Genuchten) rely on an asymptotic assumption of infinite suction at zero water content, bounded models introduce a cut-off pressure where liquid continuity breaks, ensuring the degeneracy $g(0) = 0$. For example, for functions with polynomial growth at $0$
\be{fgpq}
\lim_{u\to 0}f(u)u^{-p-1} = C_f, \qquad \lim_{u\to 0} g(u)u^{-q-1} = C_g,
\ee
with exponents $p,q$ and some positive constants $C_f, C_g$, Hypothesis~\ref{hfg1} is satisfied if and only if $p>0$ and $-1<q<2p-1$. The Rossi-Nimmo model typically implies a linear behavior $\theta \sim u$ of water content near the dry limit (corresponding to $q=0$) and a power-law growth of permeability $\kappa(\theta) \sim \theta^\gamma$ with $\gamma > 1$ (corresponding to $p = \gamma$). We observe that these physical parameters satisfy the condition $-1 < q < 2p-1$ required above, as typically $\gamma \ge 3$ for all realistic soils.

In the next section we prove the following result.

\begin{theorem}\label{tfg1}
Let Hypothesis~\ref{hfg1} hold, let $u^*>0$ be a constant, and let $u_0 \in L^\infty(\Omega)$ be given, $0 \le u_0(x) \le u^*$ a.\,e., $\nabla u_0 \in L^2(\Omega)$. Then there exists a unique solution $u\in L^r(\Omega; C[0,T])$, $1\le r < \infty$, to \eqref{1e1} with initial condition \eqref{e2} and such that $0\le u(x,t) \le u^*$ a.\,e., $f(u) \in L^2(0,T;W^{1,2}(\Omega))$, and $g(u)_t \in L^2(0,T; W^{-1,2}(\Omega))$.
\end{theorem}

In Section~\ref{sec:mois} we address the problem of bounded propagation speed and describe the shape of the moving wet region. We define the \emph{wetting front} $R^k$ in directions of the unit coordinate vectors $\bfe_k$, $k = 1, \dots, N$, as
\begin{equation}\label{wetreg}
R^k(t) \coloneqq \sup \left\{x_k : x \in {\rm supp}\,u(\,\cdot\,,t) \right\}.
\end{equation}
Note that although \eqref{wetreg} is formulated in terms of the pressure $u$, it actually determines the maximal size of the wet area (i.e., of the support of $\theta$). This is because, by Hypothesis~\ref{hfg1}, the function $g$ in \eqref{e0a} vanishes only at $0$. Therefore, we will continue to refer to the support of $u$ as the wet region.
	
The following result summarizes the main findings.

\begin{theorem}\label{3t1}
	Let the assumptions of Theorem~\ref{tfg1} be satisfied, and additionally assume that the function 
\be{ppp}
P(u) := \frac{f'(u)}{g(u)} \ \for u>0
\ee	
is uniformly bounded from above and from below.
	Let $\Omega = (-M,M)^N$ with $M>0$ sufficiently large, and let the initial condition $u_0$ have compact support in $\Omega$. Then
	\begin{enumerate}
		\item There exists a unique solution to Problem~\eqref{e0u}--\eqref{e2}, which is the same for homogeneous Dirichlet/Neumann/Robin boundary conditions, is non-negative almost everywhere, has compact support in $\Omega\times (0,T)$, and can be extended by $0$ to a global solution on $\real^N \times (0,\infty)$.
		\item The dynamics of the wetting front $R^k(t)$, $k = 1, \dots, N$, can be summarized as follows:
		\begin{enumerate}[leftmargin=5.5mm,label=(\roman*)]
			\item the \emph{lateral fronts} in directions of unit vectors $\bfe_k$, $k=1, \dots, N-1$ perpendicular to gravity propagate in time at most proportionally to $\sqrt{t}$;
			\item the \emph{downward front} in direction $-\bfe_N$ moves in time at most proportionally to $t$;
			\item the \emph{upward front} is bounded by a constant surface independent of $t$. Furthermore,
			\begin{itemize}[leftmargin=5.5mm]
				\item if $\lim_{u\to 0} P(u) > 0$, then the upward front starts moving downward after an initial time interval and achieves an asymptotically constant speed for large times;
				\item if $\lim_{u\to 0} P(u) = 0$, then the upward front reaches the upper bound asymptotically as $t\to\infty$ whenever $\int_0^{u^*} 1/P(v) \dd v = \infty$ and in finite time whenever $\int_0^{u^*} 1/P(v) \dd v < \infty$, and it does not reverse direction.
			\end{itemize}
		\end{enumerate}
	\end{enumerate}
\end{theorem}


\section{Existence of solutions}\label{sec:exis}

We first check that if a solution $u$ as in Theorem~\ref{tfg1} exists, then it is non-negative and unique.

\begin{propo}\label{1t2}
	Consider any monotone extensions $\tilde f, \tilde g: \real \to \real$ of $f,g$. Then Problem~\eqref{1e1}, \eqref{e2} for $f,g$ replaced with $\tilde f, \tilde g$ has at most one solution, and this solution, if it exists, is non-negative a.\,e. In particular, the extensions  $\tilde f, \tilde g$ are never active.
\end{propo}

\bpf{Proof}
We denote by $H$ the Heaviside function, that is, $H(s) = 0$ for $s\le 0$, $H(s) = 1$ for $s>0$, and for $\sigma > 0$ we define its regularization $H_\sigma$ by the formula
\begin{equation}\label{Heav}
	H_\sigma (s) = \left\{
	\begin{array}{ll}
		0 & \for s \le 0,\\[1mm]
		\frac{s}{\sigma} & \for s\in (0,\sigma),\\[1mm]
		1 & \for s \ge \sigma.
	\end{array}
	\right.
\end{equation}
Note that for all $u,v \in L^1(\Omega;W^{1,1}(0,T))$ we have the identity
\be{hea}
(\tilde g(v) - \tilde g(u))_t H(v-u) = \frac{\partial}{\partial t}(\tilde g(v) - \tilde g(u))^+ \ \ale \mbox{ in } \Omega.
\ee
To prove that the solution $u$ to \eqref{1e1} stays positive, we test \eqref{1e1} by $\phi = -H_\sigma (-u)$ (note that this is now an admissible test function) and integrating by parts we get
$$
\io \Big((-\tilde g(u))_t H_\sigma(-u) + \tilde f'(u) H'_\sigma(-u) |\nabla u|^2\Big)  \dd x + \io H_\sigma(-u)\frac{\partial}{\partial x_N}\tilde f'(u) \dd x = 0.
$$
The term $\tilde f'(u) H'_\sigma(-u) |\nabla u|^2$ is non-negative. Letting $\sigma \to 0$ and using the identity
$$
\io H(-u)\frac{\partial}{\partial x_N}\tilde f'(u) \dd x = \io H(-\tilde f'(u))\frac{\partial}{\partial x_N}\tilde f'(u) \dd x = 0,
$$
we thus obtain from \eqref{hea} that
$$
\io (-\tilde g(u(x,t)))^+ \dd x \le \io(-g(u(x,0)))^+ \dd x = 0,
$$
hence, $u(x,t) \ge 0$ a.\,e.

To prove uniqueness, we consider two solutions $u_1$ and $u_2$ of \eqref{1e1}, and test the difference of the two equations by $\phi = H_\sigma (f(u_1) - f(u_2))$ with $H_\sigma$ given by \eqref{Heav}. We obtain $I_1+I_2+I_3 = 0$, where
\begin{align*}
	I_1 &= \io \left(g(u_1) - g(u_2)\right)_t H_\sigma (f(u_1) - f(u_2))\dd x,\\[1mm]
	I_2 & = \io H_\sigma' (f(u_1) - f(u_2))\big| f'(u_1)\nabla u_1 - f'(u_2) \nabla u_2\big|^2\dd x,\\[1mm]
	I_3 & = \io (f'(u_1) - f'(u_2))\frac{\partial}{\partial x_N}H_\sigma (f(u_1) - f(u_2))\dd x.
\end{align*}
The integral $I_2$ is non-negative. Letting $\sigma \to 0$, using the identities $H(f(u_1) - f(u_2)) =  H(u_1 - u_2)$, and integrating in $I_3$ by parts, we obtain
$$
\io \left(g(u_1) - g(u_2)\right)_t H(u_1 - u_2)\dd x - \io H(u_1 - u_2)\frac{\partial}{\partial x_N} (f'(u_1) - f'(u_2))\dd x \le 0.
$$
From Lemma~\ref{1l1} below it follows that
$$
\io H(u_1 - u_2)\frac{\partial}{\partial x_N} (f'(u_1) - f'(u_2))\dd x = 0,
$$
and identity \eqref{hea} implies that $u_1 \le u_2$ a.\,e. Uniqueness is obtained by interchanging $u_1$ and $u_2$.
\epf

\begin{lemma}\label{1l1}
	Let $\beta: \real \to \real$ be a continuously differentiable function, and let $u_1,u_2 \in W^{1,2}(\Omega)$ be functions such that $u_1,u_2 \ge 0$ a.\,e. and $u_1 = 0$ on $\partial\Omega$. Then
	$$
	\io H(u_1 - u_2)\frac{\partial}{\partial x_N} (\beta(u_1) - \beta(u_2))\dd x = 0.
	$$
\end{lemma}

\bpf{Proof}
Let $\bar\Omega$ be the closure of $\Omega$ and let $M>0$ be such that $\bar\Omega \subset (-M,M)^N$. We extend $u_1, u_2$ on $(-M,M)^N$ in such a way that $u_1(x) = 0$ and $u_2(x) \ge 0$ for $x \in (-M,M)^N\setminus \bar\Omega$. We have
\begin{align*}
	&\io H(u_1 - u_2)\frac{\partial}{\partial x_N} (\beta(u_1) - \beta(u_2))\dd x\\[1mm] &= \int_{(-M,M)^{N-1}} \int_{-M}^M  H(u_1(x',x_N) - u_2(x',x_N))\frac{\partial}{\partial x_N} (\beta(u_1(x',x_N)) - \beta(u_2(x',x_N)))\dd x_N\dd x'.
\end{align*}
For a.\,e. $x' \in(-M,M)^{N-1}$, the functions $u_1(x',\cdot), u_2(x',\cdot)$ are absolutely continuous, and the set $A_N(x')\coloneqq\{x_N \in [-M,M]: u_1(x',x_N)\!>\! u_2(x',x_N)\}$ is open. Note that we can write $A_N(x')\! =\! \bigcup_{k=1}^\infty (a_k(x'),b_k(x'))$ with $u_1(x',a_k(x')) = u_2(x',a_k(x'))$, $u_1(x',b_k(x'))= u_2(x',b_k(x'))$, $u_1(x',x_N) > u_2(x',x_N)$ for $x_N \in (a_k(x'), b_k(x'))$. For a.\,e. $x' \in (-M,M)^{N-1}$ we then have
\begin{align*}
	&\hspace{-4mm} \int_{-M}^M  H(u_1(x',x_N) - u_2(x',x_N))\frac{\partial}{\partial x_N} (\beta(u_1(x',x_N)) - \beta(u_2(x',x_N)))\dd x_N\\[1mm]
	&= \sum_{k=1}^\infty\int_{a_k(x')}^{b_k(x')} \frac{\partial}{\partial x_N} \big(\beta(u_1(x',x_N)) - \beta(u_2(x',x_N))\big)\dd x_N\\[1mm]
	&= \sum_{k=1}^\infty \Big(\beta(u_1(x',b_k(x'))) - \beta(u_1(x',a_k(x'))) - \beta(u_2(x',b_k(x')))+\beta(u_2(x',a_k(x')))\Big)= 0,
\end{align*}
and the assertion follows.
\epf

The existence of a solution to Problem~\eqref{1e1}, \eqref{e2} with regularity as in Theorem~\ref{tfg1} will be proved in several intermediate steps. We first denote
\be{defgamma}
v = f(u),\ \  \hat g(v) = g(f^{-1}(v)),\ \  \hat f(v) = f'(f^{-1}(v)),
\ee
and rewrite \eqref{1e1} in terms of $v$ in the form
\be{1e1a}
\io \left(\hat g(v)_t \phi + \left(\nabla v + \hat f(v) \bfe_N\right) \cdot \nabla \phi\right)\dd x = 0 \quad \forall \phi \in \sobo
\ee
with initial condition
\be{iniv}
v(x,0) = v_0(x) \coloneqq f(u_0(x)).
\ee
To prove the existence of a solution to Problem~\eqref{1e1a}--\eqref{iniv}, we introduce two regularizing parameters $0<\delta <\ve < 1$ and replace \eqref{1e1a} with
\be{1e1b}
\io \left((\ve v + \gpd(v))_t \phi + \left(\nabla v + \hat f(v^+{+}\delta) \bfe_N\right) {\cdot} \nabla \phi\right)\dd x = 0 \quad \forall \phi \in \sobo,
\ee
where we denote $\gpd(v) = \hat g(v+\delta) - \hat g(\delta)$.

We first prove that if Problem~\eqref{iniv}--\eqref{1e1b} has a solution, then it is non-negative and uniformly bounded.

\begin{propo}\label{1p3}
Put $v^* = f(u^*)$ with $u^*$ from Theorem~\ref{tfg1}, and let $v= \vved$ be a solution to Problem~\eqref{iniv}--\eqref{1e1b} such that $v_t \in L^2(\Omega\times (0,T))$ and  $\nabla v \in L^\infty(0,T; L^2(\Omega))$. Then
\be{vvK}
0\le v(x,t) \le v^* \ \ale
\ee
\end{propo}

\bpf{Proof}
Let us first prove that $v = \vved \ge 0$ a.\,e.  We can test \eqref{1e1b} directly by $\phi = -v^-$, which yields that
$$
\frac{\dd}{\dd t}\io \mathcal{E}_\ve((\vved)^-)\dd x +  \io |\nabla (\vved)^-|^2 \dd x = \hat f(\delta) \io \frac{\partial (\vved)^-}{\partial x_N} \dd x = 0
$$
with $\mathcal{E}_\ve(v) = \frac{\ve}{2}v^2 + \int_0^v y \gpd'(-y)\dd y$, and the assertion follows from the positivity of the initial condition.
\\
To obtain the upper bound, we test \eqref{1e1b} by $\phi = H_\sigma (v - v^*)$ with $H_\sigma$ given by \eqref{Heav} and proceed as in the proof of Proposition~\ref{1t2}. More specifically, we have
$$
\nabla v \cdot \nabla H_\sigma (v - v^*)= |\nabla v|^2 H'_\sigma (v - v^*) \ge 0 \ \ale,
$$
$$
\io \hat f(v+\delta)\frac{\partial}{\partial x_N}H_\sigma (v - v^*)\dd x = -\io H_\sigma (v - v^*) \frac{\partial }{\partial x_N}\hat f(v+\delta)\dd x.
$$
Using Lemma~\ref{1l1} and arguing as in the proof of Proposition~\ref{1t2}, we obtain the identity
$$
\lim_{\sigma \to 0} \io \hat f(v+\delta)\frac{\partial}{\partial x_N}H_\sigma (v - v^*)\dd x = 0.
$$
We conclude that, for all $t\in (0,T)$,
$$
\io \big(\ve (v-v^*) + (\gpd(v) - \gpd(v^*))\big)_t \,H(v - v^*)\dd x \le 0.
$$
By virtue of \eqref{hea} we have
$$
\io \big(\ve (v{-}v^*) + (\gpd(v) {-} \gpd(v^*))\big)_t H(v {-} v^*)\dd x = \frac{\dd}{\dd t}\io \left(\ve (v{-}v^*)^+ + (\gpd(v) {-} \gpd(v^*))^+\right)\dd x \le 0,
$$
and from \eqref{vvK} we obtain $v(x,t) \le v^*$ a.\,e., which we wanted to prove.
\epf

We now check that Problem~\eqref{iniv}--\eqref{1e1b} has a solution $v=\vved$ with regularity independent of $\ve$ and $\delta$.

\begin{propo}\label{2l1}
Under the hypotheses of Theorem~\ref{tfg1}, Problem~\eqref{iniv}--\eqref{1e1b} has for every $\delta>0$ and $\ve > 0$ a solution $v = \vved$ such that $v_t \in L^2(\Omega\times (0,T))$ and  $\nabla v \in L^\infty(0,T; L^2(\Omega))$.
\end{propo}

\bpf{Proof}
For the purpose of Proposition~\ref{2l1}, we can assume that both $\hat f$ and $\gpd$ have linear growth at infinity. Indeed, by Proposition~\ref{1p3} the values of $v$ are in the interval $[0,v^*]$, and thus, the values of the nonlinearities outside this interval are not relevant.

Let $\{w_k: k\in \nat\} \subset L^2(\Omega)$ be the complete orthonormal system of eigenfunctions of the negative Laplace operator with the boundary condition \eqref{e2}, that is,
\be{ortho}
-\Delta w_k = \mu_k w_k, \quad w_k = 0 \ \mbox{ on } \ \partial\Omega,
\ee
and for $n \in \nat$ and given functions $v_1, \dots, v_n \in C^1[0,T]$ put $v\on(x,t) = \sum_{k=1}^n v_k(t) w_k(x)$. We consider a finite-dimensional Galerkin approximation of \eqref{1e1b} in the form
\be{2e1g}
\ve \dot v_k + \mu_k v_k +\io \left(\big(\gpd(v\on)\big)_t w_k + \left(\hat f((v\on)^+{+}\delta)\right)  \frac{\partial w_k}{\partial x_N}\right)\dd x  = 0
\ee
for $k=1, \dots, n$ with initial conditions $v_k(0) = \io v(x,0) w_k(x)\dd x$. This is a simple ODE system which admits at least a local solution $\{v_1, \dots, v_n\} \subset C^1[0,T_n]$, $T_n \le T$.

We now derive some estimates which will enable us to conclude that $T_n = T$ and pass to the limit as $n \to \infty$. The first step is the energy estimate, which is obtained by multiplying \eqref{2e1g} by $v_k$ and summing over $k=1, \dots, n$. The terms under the time derivative are estimated by the initial condition. The last term under the integral sign in \eqref{2e1g} vanishes because of the zero Dirichlet boundary condition, and we obtain
\be{est1}
\int_0^{T_n}\io |\nabla v\on|^2 \dd x \le C
\ee
with a constant $C>0$ independent of $n$, $\delta$, and $\ve$.

Next, we multiply \eqref{2e1g} by $\dot v_k$, sum up over $k=1, \dots, n$, and integrate by parts to get
\be{von}
\io \ve |v\on_t|^2\dd x + \frac{1}{2}\frac{\dd}{\dd t}\io |\nabla v\on|^2\dd x \le \io v\on_t \hat f'((v\on)^+ +\delta)\frac{\partial v\on}{\partial x_N}\dd x.
\ee
The right-hand side of \eqref{von} can be estimated from above by the term
\be{hold}
C_\delta\left(\io |v\on_t|^2 \dd x\right)^{1/2} \left(\io |\nabla v\on|^2\dd x\right)^{1/2}
\ee
with a constant $C_\delta>0$ independent of $n$ and $\ve$. Hence, by Young's inequality and \eqref{est1}, there exists a constant $C_{\ve,\delta}{>} 0$ independent of $n$ and depending possibly on $\ve$ and $\delta$ such that
\be{ested}
\int_0^{T_n}\io |v\on_t|^2 \dd x\dd t + \sup_{t\in (0,T_n)}\io |\nabla v\on|^2 \dd x \le C_{\ve,\delta},
\ee
where we have also used the assumptions on $u_0$. Therefore, $T_n=T$ and, up to subsequences, $v\on$ converges to some $v = \vved$ strongly in $L^2(\Omega \times (0,T))$, $v\on_t$ converges to $\vved_t$ weakly in $L^2(\Omega \times (0,T))$, and $\nabla v\on$ converges to $\nabla \vved$ weakly in $L^2(\Omega \times (0,T))$ and weakly-star in $L^\infty(0,T;L^2(\Omega))$ as $n \to \infty$. We conclude that $v =\vved$ is a solution to the PDE
\be{1e1c}
\io \left((\ve v + \gpd(v))_t \phi + \left(\nabla v + \hat f(v{+}\delta) \bfe_N\right) {\cdot} \nabla \phi\right)\dd x = 0 \quad \forall \phi \in \sobo,
\ee
which completes the proof.\epf

We now show that the solution $v$ to Problem~\eqref{1e1a}--\eqref{iniv} is obtained by letting $\delta$ and $\ve$ tend to $0$.

\begin{theorem}\label{t0}
Under the hypotheses of Theorem~\ref{tfg1}, Problem~\eqref{1e1a}--\eqref{iniv} has a unique solution $v\in L^p(\Omega; C[0,T])$  for every $1\le p < \infty$ such that $0 \le v(x,t) \le v^*$ a.\,e.\ with $v^*$ from Proposition~\ref{1p3}, $\nabla v \in L^\infty(0,T;L^2(\Omega))$, and $\hat g (v)_t \in L^2(0,T;W^{-1,2}(\Omega))$.
\end{theorem}

\bpf{Proof}
Uniqueness follows from Proposition~\ref{1t2}. The existence proof is carried out by passing to the limit in \eqref{1e1b} as $\delta \to 0$ and $\ve \to 0$. To this end, we need estimates independent of $\ve$ and $\delta$. For simplicity, we write again $v$ instead of $\vved$ and introduce the functions
\be{dh}
h(u) = \int_0^u \frac{|f''(z)|^2}{g'(z)}\dd z \ \for u\ge 0, \quad \hat h(v) = h(f^{-1}(v)) \ \for v \ge 0.
\ee
These are indeed well defined functions by Hypothesis~\ref{hfg1}.
We now test \eqref{1e1b} by $\phi = \hat h_\delta(v) \coloneqq \hat h(v+\delta) - \hat h(\delta)$, which is an admissible test function, and we get the identity
$$
\frac{\dd}{\dd t} \io G_{\ve,\delta}[v]\dd x +  \io |\nabla v|^2\hat h'(v+\delta)\dd x = \io \frac{\partial}{\partial x_N} \hat F(v)\dd x = 0,
$$
where
$$
\hat F(v) = \int_0^{v} \hat f(z+\delta) \hat h'_\delta(z)\dd z, \quad
G_{\ve,\delta}[v] = \int_0^v (\ve + \hat g_\delta'(z+\delta))\hat h_\delta(z)\dd z.
$$
The function $G_{\ve,\delta}[v](x,0)$ is bounded independently of $\ve$ and $\delta$ and we get
\be{est2b}
\int_0^T\io |\nabla \vved|^2\hat h'(\vved+\delta) \dd x\dd t \le C
\ee
with a constant $C>0$.

To estimate the dependence of $v_t$ on $\ve$ and $\delta$, we denote for $k,n \in \nat$
$$
\hat v_k(t) =\io v(x,t) w_k(x) \dd x, \qquad \hat v\on (x,t) = \sum_{k=1}^n \hat v_k(t) w_k(x)
$$
with $w_k$ as in \eqref{ortho}, test \eqref{1e1b} with $\phi = \hat v\on_t$ and integrate by parts. Letting $n \to \infty$ we get for any $t \in (0,T)$ a counterpart of \eqref{von}, namely,
\be{voe}
\int_0^t\io \big(\ve + \hat g'(v{+}\delta)\big) |v_t|^2\dd x\dd\tau + c\!\io \!|\nabla v|^2(x,t)\dd x \le C\!\left(\!1 {+}\!\!\int_0^t\io |v_t| \hat f'(v {+} \delta)\left|\frac{\partial v}{\partial x_N}\right|\dd x\dd\tau\!\right)
\ee
with positive constants $c,C$ depending only on the initial data. We estimate the right-hand side of \eqref{voe} using H\"older's inequality as
\be{hoel}
\io |v_t| \hat f'(v {+} \delta)\left|\frac{\partial v}{\partial x_N}\right|\dd x \le \left(\io \hat g'(v {+} \delta)|v_t|^2\dd x\right)^{1/2} \left(\io \frac{|\hat f'(v{+}\delta)|^2}{\hat g'(v{+}\delta)}|\nabla v|^2\dd x\right)^{1/2}.
\ee
By \eqref{defgamma} and \eqref{dh} we have
$$
\frac{|\hat f'(v+\delta)|^2}{\hat g'(v+\delta)} = \hat h'(v+\delta),
$$
and we obtain from \eqref{voe} and \eqref{est2b} for solutions $v = \vved$ of \eqref{1e1b} bounds independent of $\ve$ and $\delta$ of the form
\be{eed}
\int_0^t\io \big(\ve + \hat g'(\vved)\big) |\vved_t|^2\dd x\dd\tau + c\io |\nabla \vved|^2(x,t)\dd x \le C.
\ee

Keeping $\ve$ constant for the moment, we can now let $\delta \to 0$ and conclude that $v^\ve = \lim_{\delta \to 0}\vved$ is a solution to the problem
\be{ev1}
\io \left(\big(\ve v^\ve + \hat g(v^\ve)\big)_t\phi + \left(\nabla v^\ve + \hat f(v^\ve) \bfe_N\right)\cdot \nabla\phi\right)\dd x = 0 \quad \forall \phi\in W^{1,2}_0(\Omega)
\ee
and satisfies estimates \eqref{eed} and \eqref{vvK} independently of $\varepsilon$. These estimates are preserved in the limit as $\delta \to 0$, that is,
\be{eee}
\int_0^t\io \big(\ve + \hat g'(v^\ve)\big) |v^\ve_t|^2\dd x\dd\tau + c\io |\nabla v^\ve|^2(x,t)\dd x \le C, \quad 0\le v^\ve(x,t) \le v^* \ \ale
\ee
We distinguish the cases
\begin{itemize}
\item[(1)] $\displaystyle{\lim_{v \to 0+} \hat g'(v) = 0}$,
\item[(2)] $\displaystyle{\hat g'(v) \ge m}$ for some $m > 0$ and for all $v \in (0,v^*)$.
\end{itemize}
In case (1) we put $w^\ve = \int_0^{v^\ve} \sqrt{\hat g'(z)} \dd z$ and obtain for all $t\in (0,T)$ the inequality
\be{eew}
\int_0^t\io \left((\ve + \hat g'(v^\ve))|v^\ve_t|^2+|w^\ve_t|^2\right)\dd x\dd\tau + \io \left(|\nabla v^\ve|^2+|\nabla w^\ve|^2\right) (x,t)\dd x \le C
\ee
with a constant $C>0$ independent of $\ve$ and $t$. By anisotropic embedding (see \cite{bin}), and passing to a subsequence if necessary, we get $w^\ve \to w$ strongly in $L^{p}(\Omega; C[0,T])$ for every $p\ge 1$. Putting $\Phi(v) \coloneqq \int_0^{v} \sqrt{\hat g'(z)} \dd z$, we see that the function $\Phi$ is continuous and increasing in $[0,v^*]$, $\Phi(0)=0$, and $\Phi(v)>0$ for $v>0$. Hence, $\Phi^{-1} : [0,\Phi(v^*)] \to [0,v^*]$ is continuous. Therefore, by dominated convergence, $v^\ve \to v$ strongly in $L^{p}(\Omega \times (0,T))$ for every $p\ge 1$. From \eqref{eew} we also obtain $w^\ve_t \to w_t$ and $\hat g'(v^\ve) v^\ve_t \to \hat g'(v) v_t$ weakly in $L^2(\Omega\times (0,T))$, $\nabla w^\ve \to \nabla w$ and $\nabla v^\ve \to \nabla v$ weakly* in $L^\infty(0,T; L^2(\Omega))$, $\ve v^\ve_t \to 0$ strongly in $L^2(\Omega\times (0,T))$, and we conclude that $v$ is a solution to Problem~\eqref{1e1a}--\eqref{iniv}.

In case (2) we have
\be{eewv}
\int_0^t\io \left((\ve + \hat g'(v^\ve) + m)|v^\ve_t|^2\right)\dd x\dd\tau + \io |\nabla v^\ve|^2 (x,t)\dd x \le C
\ee
with a constant $C>0$ independent of $\ve$ and $t$. We find again a subsequence such that $v^\ve \to v$ strongly in $L^{p}(\Omega; C[0,T])$ for every $p\ge 1$, $v^\ve_t \to v_t$ weakly in $L^2(\Omega\times (0,T))$, $\nabla v^\ve \to \nabla v$ weakly* in $L^\infty(0,T; L^2(\Omega))$, $\ve v^\ve_t \to 0$ strongly in $L^2(\Omega\times (0,T))$. Again by comparison in \eqref{ev1}, we see that the terms $\hat g(v^\ve)_t$ are bounded in $L^2(0,T; W^{-1,2}(\Omega))$ independently of $\ve$ and $\hat g(v^\ve) \to \hat g(v)$ strongly in $L^{p}(\Omega; C[0,T])$ for every $p\ge 1$. Passing to the limit in \eqref{ev1}, we complete the proof of Theorem~\ref{t0}.
\epf


\section{Moisture propagation speed}\label{sec:mois}

In this section we prove Theorem~\ref{3t1}. We fix a constant $M>0$ and consider the solution $u$ to Problem~\eqref{1e1}, \eqref{e2} on
\be{ome}
\Omega = (-M,M)^N
\ee
under the assumption that, for some $R_0 \in (0,M)$ and $u^*>0$, we have
\be{r0}
0 \le u(x,0) = u_0(x) \le u^* \ \for\ale \ x \in \Omega, \quad u_0(x) = 0 \ \for \ |x| > R_0.
\ee
In other words, it is assumed that no liquid is initially present in the body outside the ball centered at $0$ of radius $R_0$.

As an easy result which will be refined later on, we first check that the moisture front can never reach the upper boundary $x_N = M$ of $\Omega$ if $M$ is sufficiently large.

\begin{propo}\label{xp3}
Let \eqref{ome} and \eqref{r0} hold, and let $u$ be the solution to Problem~\eqref{1e1}, \eqref{e2}. Then we have
$$
u(x,t) \le \big(R_N^\flat - x_N\big)^+\ \ \ale, \quad R_N^\flat \coloneqq u^* + R_0.
$$
\end{propo}

\bpf{Proof}
For $x \in \Omega$ put $U(x) = (R_N^\flat -x_N)^+$. We have
$$
f'(U)(\nabla U + \bfe_N) = f'((R_N^\flat -x_N)^+)\big(1 - H(R_N^\flat -x_N)\big)\bfe_N = 0
$$
for a.\,e. $x \in \Omega$. Hence, in particular,
\be{1e0}
\io \left(g(u)_t \phi + \Big(f'(u) (\nabla u + \bfe_N) - f'(U)(\nabla U + \bfe_N)\Big)\cdot \nabla \phi\right)\dd x = 0
\ee
for all $\phi\in W^{1,2}_0(\Omega)$. Choosing $\phi = H_\sigma (f(u) - f(U))$ with $H_\sigma$ given by \eqref{Heav} for $\sigma>0$, letting $\sigma \to 0$, and arguing as in the proof of Proposition~\ref{1t2} we get for all $t\in (0,T)$ that
$$
\io \big(g(u(x,t)) - g(U(x))\big)_t \,H(u(x,t) - U(x))\dd x \le 0,
$$
and from \eqref{hea} we obtain
$$
\frac{\dd}{\dd t}\io\left(\big(g(u(x,t)) - g(U(x))\big)^+ \right)\dd x \le 0.
$$
The choice of $R_N^\flat$ implies that $u_0(x) \le U(x)$ a.\,e., hence,
$$
\io\left(\big(g(u(x,t)) - g(U(x))\big)^+ \right)\dd x \le 0,
$$
which we wanted to prove.
\epf

\subsection{Traveling wave solutions}

To study the moisture front propagation in other directions, we construct dominant traveling wave solutions $u_\omega$ of \eqref{1e1} in each direction
\be{eome}
\bfe_\omega = \som\, \bfe_\perp - \com\, \bfe_N
\ee
for an arbitrarily fixed $\omega \in [-\pi/2,3\pi/2)$ and a unit vector $\bfe_\perp$ orthogonal to $\bfe_N$, see Figure~\ref{f1}. For $x_\perp, x_N \in \real$, we denote $x = x_\perp \bfe_\perp + x_N \bfe_N$ and assume $u_\omega$ in the form
\be{uom}
u_\omega(x,t) = U_{c,\omega}(ct + R_\omega - \scal{x,\bfe_\omega}) = U_{c,\omega}(ct + R_\omega - x_\perp\,\som + x_N\,\com)
\ee
with a differentiable function $U_{c,\omega}: \real \to \real^+$ such that $U_{c,\omega}(z) = 0$ for $z\le 0$, and with a constant shift $R_\omega$ which will be specified below in \eqref{rom}. We have for $z=ct + R_\omega - \scal{x,\bfe_\omega}$ that
$$
\nabla u_\omega(x,t) = -\frac{\dd}{\dd z}U_{c,\omega}(z)\bfe_\omega, \quad \dive \big(f'(u_\omega(x,t))\nabla u_\omega(x,t)\big) = \frac{\dd}{\dd z}\left(f'(U_{c,\omega})\frac{\dd}{\dd z}U_{c,\omega}(z)\right)
$$
hence, $u_\omega$ is a solution to \eqref{e0u} if and only if for all $z\ge 0$ we have
\be{odeom}
c g(U_{c,\omega}(z)) =  f'(U_{c,\omega}(z))(U'_{c,\omega}(z) + \com).
\ee
Let $P(u)$ be given by \eqref{ppp} and let us assume that there exist constants $0\le P^\flat\le P_0 \le P^\sharp$ such that
\be{pp0}
\lim_{u\to 0+}P(u) = P(0) = P_0, \quad P^\flat \le  P(u) \le P^\sharp \quad \forall u\in (0,u^*].
\ee
We define $c^*_\omega = P^\sharp\com$ if $\com > 0$, $c^*_\omega = -P^\flat |\com|$ if $\com < 0$, and $c^*_0 = 0$, and put
\be{defq}
Q_{c,\omega}(u) \coloneqq \int_0^u \frac{P(v)}{c- P(v)\com} \dd v \ \for u \in [0,u^*] \mathrm{\ and \ } c>c^*_\omega.
\ee
Note that in case \eqref{fgpq}, condition \eqref{pp0} is satisfied if and only if $p\ge q+1$. This certainly holds for all realistic soils, see the discussion following Hypothesis~\ref{hfg1}.

With $Q_{c,\omega}(u)$ as in \eqref{defq} we thus have
\be{pzom}
\frac{\dd}{\dd z}Q_{c,\omega}(U_{c,\omega}(z)) = 1,
\ee
or, in other terms,
\be{qzom}
U_{c,\omega}(z) = Q_{c,\omega}^{-1}(z)\ \for z>0.
\ee
This proves that, if condition \eqref{pp0} holds, traveling wave solutions as in \eqref{uom} can be constructed via the identity \eqref{qzom}. We now need to prove that they are dominant. To guarantee that $u_\omega(x,0) \ge u_0(x)$ a.\,e., it is enough to check that the implication
\be{impl}
|x| \le R_0 \ \Longrightarrow \ R_\omega - \scal{x,\bfe_\omega} \ge Q_{c,\omega}(u^*)
\ee
holds, which is certainly the case if we put
\be{rom}
R_\omega = R_0 + Q_{c,\omega}(u^*).
\ee
We now proceed as in the proof of Proposition~\ref{1t2}, test the difference
$$
\io \big(g(u) - g(u_\omega)\big)_t\phi \dd x + \io \big(f'(u)(\nabla u + \bfe_N) - f'(u_\omega)(\nabla u_\omega+\bfe_N)\big)\cdot \nabla \phi \dd x = 0\qquad \forall\phi\in W^{1,2}_0(\Omega)
$$
by $\phi = H_\sigma(f(u)- f(u_\omega))$, and let $\sigma$ tend to $0$ to get the estimate
\be{bom}
u(x,t) \le u_\omega(x,t) \ \ale
\ee


\subsection{Shape of the moving wet region}\label{shape}

The coordinate $\rho$ of the moving moisture front $x = \rho \bfe_\omega$ in direction $\bfe_\omega$ is characterized by the equation $ct + R_\omega - \scal{x, \bfe_\omega} = 0$, that is, in view of \eqref{rom} and \eqref{defq},
\be{eom}
ct + R_0 + \int_0^{u^*} \frac{P(v)}{c- P(v)\com} \dd v = \rho, \quad c>c^*_\omega.
\ee


\subsubsection{Lower half-space}

Assume first that $\com \ge 0$, that is, $\omega \in [-\pi/2,\pi/2]$. By virtue of \eqref{defq}, there are no traveling waves in direction $\bfe_\omega$ with propagation speed lower than $P^\sharp\com$. The \emph{envelope} of the straight lines given by \eqref{eom} is described by the differential equation
\be{eot}
t\dot R^\omega(t) + R_0 + \int_0^{u^*} \frac{P(v)}{\dot R^\omega(t)-P(v)\com} \dd v = R^\omega(t), \quad R^\omega(0) = R_0.
\ee
Note that \eqref{eot} is the so-called Clairaut equation (see \cite[Chapter~XIII, Sec.~13]{pisk}) for the unknown $y(t) = R^\omega(t) - R_0$ of the form
\be{eot1a}
t\dot y(t) + \psi(\dot y(t)) = y(t), \quad y(0) = 0,
\ee
with a decreasing convex function $\psi:(\bar y, \infty) \to \real^+$ of class $C^\infty$ such that $\psi(\infty) = 0$, $0<\psi(\bar y+) \leq \infty$ with $\bar y=c^*_\omega=P^\sharp\com$. A solution $y:[0,\infty) \to \real$ of \eqref{eot1a} is therefore necessarily positive, increasing, $y(t) > 0$ and $\dot y(t) \ge \bar y$ for $t>0$, and $\lim_{t\to 0+} \dot y(t) = \infty$.

Differentiating \eqref{eot1a} with respect to $t$ yields
\be{eot2a}
\ddot y(t)\left(t + \psi'(\dot y(t)) \right) = 0.
\ee
Hence, a solution to \eqref{eot1a} is constructed by concatenating the so-called general solution (a family of straight lines $\ddot y(t) = 0$) with the singular solution defined by the ODE
\be{eot2b}
t + \psi'(\dot y(t)) = 0.
\ee
Since $\dot y$ is unbounded near $t=0$, also $\ddot y$ is unbounded in a right neighborhood of $0$, and the system must initially follow the singular branch. Let us assume that there exist $0<a<b<\infty$ such that $\ddot y(b)=0$, $\ddot y(t) \ne 0$ for $t\in (a,b)$. Differentiating the singular condition \eqref{eot2b} gives
\be{eot3}
1+ \ddot y(t) \psi''(\dot y(t)) = 0
\ee
for $t \in (a,b)$, and we conclude that $\dot y(b-) = \bar y$ and $\ddot y(t) < 0$ for all $t\in (a,b)$. In particular, there are no intermediate points $\tau \in (0,b)$ with $\ddot y(\tau) = 0$. Hence, $\ddot y(t) < 0$ for all $t\in (0,b)$, the singular phase ends at time $b$, and the solution extends linearly for $t \ge b$ as $y(t) = y(b) + (t-b)\bar y$. The solution is always of class $C^1(0, \infty)$. However, at the transition time $t=b$, a discontinuity in the second derivative occurs unless $\psi''(\bar y+) = \infty$. Specifically, $\lim_{t\to b-} \ddot y(t) = -1/\psi''(\bar y+)$ while $\ddot y(t) = 0$ for $t > b$. If such an interval $(a,b)$ does not exist, then the singular phase is active on the whole interval $(0,\infty)$.

In terms of the original problem \eqref{eot}, the counterpart of the singular condition \eqref{eot2b} reads
\be{eot2}
\int_0^{u^*} \frac{P(v)}{(\dot R^\omega(t)-P(v)\com)^2} \dd v = t.
\ee
Putting $Q^* = \int_0^{u^*}P(v)\dd v$ we have
\be{lim0}
\lim_{t\to 0+} \dot R^\omega(t) = +\infty, \qquad \lim_{t\to 0+} \sqrt{t} \dot R^\omega(t) = \sqrt{Q^*}.
\ee
The case $\com = 0$ is straightforward and yields
\be{rom0}
R^{\pm\pi/2}(t) = R_0 + 2\sqrt{Q^*\,t} \quad \mbox{for all } \ t>0.
\ee
In other words, for small times in all directions, and for large times in directions perpendicular to gravity, the distance of the moisture front from the origin increases proportionally to $\sqrt{t}$ at most.

Let us now assume that $\com > 0$, and define the critical time $T^\sharp$ as
\be{tsha}
T^\sharp \coloneqq \frac{1}{\cos^2\omega}\int_0^{u^*} \frac{P(v)}{(P^\sharp-P(v))^2} \dd v.
\ee
According to the above discussion, one of the following two cases occurs:
\begin{itemize}
	\item[{\rm a.}] $T^\sharp = \infty$. Then Eq.~\eqref{eot2} determines the solution $R^\omega(t)$ for all $t \in (0,\infty)$, and $\lim_{t\to \infty} \dot R^\omega(t) = P^\sharp\com$;
	\item[{\rm b.}] $T^\sharp < \infty$. Then Eq.~\eqref{eot2} determines the solution $R^\omega(t)$ only for $t \in (0,T^\sharp)$ with $\dot R^\omega(T^\sharp-) = P^\sharp\com$, and on $[T^\sharp,\infty)$ it can be extended by the linear function $R^\omega(t) = R^\omega(T^\sharp) + (t-T^\sharp)P^\sharp\com$.
\end{itemize}
In other words, the downward wetting front reaches the speed $P^\sharp\com$ asymptotically as $t\to\infty$ in Case~a, and in finite time in Case~b.

\begin{figure}[htb]
\begin{center}
\includegraphics[height=5cm]{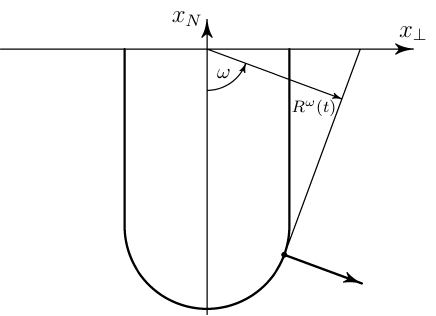}
\caption{Intersection of the downward wetting front with the plane $E(\bfe_\perp,\bfe_N)$.}
\label{f1}
\end{center}
\end{figure}

Putting
\be{romt2}
\Gamma_\omega(s) = \int_0^{u^*} \frac{P(v)}{\bigl(s - P(v)\com\bigr)^2} \dd v \quad \for \ s \in (P^\sharp\com, \infty)
\ee
and using \eqref{eot}, we can represent the solution $R^\omega(t)$ to \eqref{eot2} by the formula
\begin{align} \label{romt}
R^\omega(t) &= R_0 + \int_0^t \Gamma_\omega^{-1}(\tau) \dd\tau \quad \for \ 0\leq t < T^\sharp,\\[1mm] \label{yan}
R^\omega(t) &= R_0 + \frac{1}{\com}\int_0^{u^*} \frac{P(v)}{P^\sharp - P(v)} \dd v + t P^\sharp \com \quad \for \ t\geq T^\sharp.
\end{align}
We can rewrite \eqref{romt2} using the identities
$$
\Gamma_\omega(s) = \frac{1}{\cos^2\omega}\Gamma_0\left(\frac{s}{\com}\right), \quad \Gamma_0(z) = \int_0^{u^*} \frac{P(v)}{(z-P(v))^2} \dd v \quad \for z > P^\sharp.
$$
From \eqref{eot2} it follows that
$$
\Gamma_0\left(\frac{\dot R^\omega(t)}{\com}\right) = t \cos^2\omega \quad \for \ 0\leq t < T^\sharp,
$$
hence,
\be{romdot}
\dot R^\omega(t) = \com\,\Gamma_0^{-1}(t \cos^2\omega), \quad R^\omega(t) = R_0 + \frac{1}{\com}\int_0^{t\cos^2\omega} \Gamma_0^{-1}(s)\dd s, \quad \for \ 0\leq t < T^\sharp.
\ee
Since $\Gamma_0^{-1}(s) \ge P^\sharp$ for all $0<s < T^\sharp$, we immediately get from \eqref{romdot} the inequality
\be{rom2}
 R^\omega(t) \ge R_0 + t\,P^\sharp\com \quad\text{for } 0\leq t < T^\sharp.
\ee
On the other hand, we have for all $z\in (P^\sharp, \infty)$ that
$$
\frac{Q^*}{z^2} \le \Gamma_0(z) \le \frac{Q^*}{(z-P^\sharp)^2},
$$
hence,
$$
\frac{Q^*}{(\Gamma_0^{-1}(\tau))^2} \le \tau \le \frac{Q^*}{(\Gamma_0^{-1}(\tau)-P^\sharp)^2} \quad \forall \ 0<\tau < T^\sharp\cos^2\omega,
$$
that is,
\be{g0t}
\sqrt{\frac{Q^*}{\tau}} \le \Gamma_0^{-1}(\tau) \le P^\sharp + \sqrt{\frac{Q^*}{\tau}} \quad \forall \ 0<\tau < T^\sharp\cos^2\omega.
\ee
From \eqref{romdot}--\eqref{g0t} we thus obtain the estimate
\be{esro}
R_0 + \max\big\{t\,P^\sharp\com,\,2\sqrt{Q^* t}\big\}  \le R^\omega(t) \le R_0 + t\,P^\sharp\com + 2\sqrt{Q^* t} \quad\text{for } 0\leq t < T^\sharp.
\ee
Combining \eqref{esro} with \eqref{yan} we thus get, in addition to \eqref{lim0}, in both cases $T^\sharp = \infty$ and $T^\sharp < \infty$
the exact convergence rate
\be{limc}
\lim_{t\to \infty}\frac{R^\omega(t)}{t} = P^\sharp\com.
\ee
Recall that at time $t$ the traveling wave $u_\omega(x,t)$ given by \eqref{uom} vanishes for $\scal{x,\bfe_\omega} \ge \rho$ with $\rho$ given by formula \eqref{eom}. Hence, by virtue of \eqref{bom}, \eqref{yan}, and \eqref{esro}, we have the implication
\be{bom2}
\scal{x,\bfe_\omega} > R^\omega_\sharp(t)\coloneqq R_0 + t\,P^\sharp\com + 2\sqrt{Q^* t} \ \Longrightarrow \ u(x,t) = 0.
\ee
For each fixed time $t$, the wet region $W(t) = \{x \in (-M, M)^N: u(x,t) > 0\}$ is therefore contained in the intersection of the half-spaces
\be{bom3}
W(t) \subset W_+(t)\coloneqq \bigcap_{\omega\in (-\pi/2,\pi/2)} V_\omega(t), \quad V_\omega(t) = \{x \in \real^N: \scal{x,\bfe_\omega} \le R^\omega_\sharp(t)\}
\ee
over all $\bfe_\perp$ perpendicular to $\bfe_N$. For each particular choice of $\bfe_\perp$, we can describe the intersection of $W_+(t)$ with the plane $E(\bfe_\perp,\bfe_N)$ generated by $\bfe_\perp$ and $\bfe_N$ explicitly. Indeed, with $x_N = \scal{x,\bfe_N}$ and $\scal{x_\perp, \bfe_N} = 0$, the condition
$$
x_\perp \som - x_N \com \le R_0 + t\,P^\sharp\com + 2\sqrt{Q^* t} \quad \forall \omega\in (-\pi/2,\pi/2)
$$
in \eqref{bom3} can be rewritten as
$$
x_\perp \som - (x_N+t\,P^\sharp) \com \le R_0 +2\sqrt{Q^* t} \quad \forall \omega\in (-\pi/2,\pi/2),
$$
which is in turn equivalent to (see Figure~\ref{f1})
\begin{align}\label{bound2}
\mbox{either }\  &x_N \ge -t\,P^\sharp \ \mbox{ and } \ |x_\perp|\le R_0 +2\sqrt{Q^* t},\\[2mm] \label{bound1}
\mbox{or }\  &x_N < -t\,P^\sharp \ \mbox{ and } \ |x_\perp|^2 + \left|x_N+ t\,P^\sharp\right|^2\le (R_0 +2\sqrt{Q^* t})^2.
\end{align}


\subsubsection{Upper half-space}

The case of $\cos \omega<0$, that is, $\omega \in (\pi/2, 3\pi/2)$, is more delicate.
Repeating the argument leading to \eqref{eot} we get for $R^\omega(t)$ the differential equation
\be{eotp}
t\dot R^\omega(t) + R_0 + \int_0^{u^*} \frac{P(v)}{\dot R^\omega(t)+P(v)|\com|} \dd v = R^\omega(t), \quad R^\omega(0) = R_0,
\ee
and also
\be{eot4}
\int_0^{u^*} \frac{P(v)}{(\dot R^\omega(t)+P(v)|\com|)^2} \dd v = t,
\ee
with the same asymptotic behavior as in \eqref{lim0} for $t\to 0$.

Similarly as in \eqref{tsha}, put
\be{tfla}
T^\flat \coloneqq \frac{1}{\cos^2\omega}\int_0^{u^*} \frac{P(v)}{(P(v) - P^\flat)^2} \dd v.
\ee
Then in both cases $T^\flat = \infty$ and $T^\flat < \infty$ we have
\be{limi}
\lim_{t\to T^\flat} \dot R^\omega(t) = P^\flat \com = -P^\flat|\com|,
\ee
and the upward wetting front reaches the speed \eqref{limi} either asymptotically or in finite time.
We now proceed as in \eqref{romt2}--\eqref{yan} and put
\be{romt2a}
\hat\Gamma_\omega(s) = \int_0^{u^*} \frac{P(v)}{(s + P(v)|\com|)^2} \dd v \quad \for \ s> -P^\flat|\com|.
\ee
Then
\begin{align} \label{romta}
R^\omega(t) &= R_0 + \int_0^t \hat\Gamma_\omega^{-1}(\tau) \dd\tau \quad \for \ 0\leq t < T^\flat,\\[1mm] \label{yan3}
R^\omega(t) &= R_0+ \frac{1}{|\com|}\int_0^{u^*} \frac{P(v)}{P(v)- P^\flat} \dd v - t P^\flat |\com| \quad \for \ t\geq T^\flat.
\end{align}
We can rewrite \eqref{romt2a} using the identities
\be{gamo}
\hat\Gamma_\omega(s) = \frac{1}{\cos^2\omega}\hat\Gamma_0\left(\frac{s}{|\com|}\right), \quad \hat\Gamma_0(z) = \int_0^{u^*} \frac{P(v)}{(z+P(v))^2} \dd v \quad \for \ z>-P^\flat.
\ee
From \eqref{eot4} it follows that
\be{romdo}
\hat\Gamma_0\left(\frac{\dot R^\omega(t)}{|\com|}\right) = t \cos^2\omega \quad \for \ 0\leq t < T^\flat,
\ee
hence,
\be{romdota}
\dot R^\omega(t) = |\com|\,\hat\Gamma_0^{-1}(t \cos^2\omega), \quad R^\omega(t) = R_0 + \frac{1}{|\com|}\int_0^{t\cos^2\omega} \hat\Gamma_0^{-1}(s)\dd s, \quad \for \ 0\leq t < T^\flat.
\ee
The function $R^\omega(t)$ increases only as long as $\hat\Gamma_0^{-1}(t \cos^2\omega) \ge 0$, that is, for $t \in (0,T^0)$, where $T^0 \coloneqq \hat\Gamma(0)/\cos^2\omega \le T^\flat$. As a refinement of Proposition~\ref{xp3}, we get the inequality
\begin{align}\nonumber
R^\omega(t) &\le R_0 + \frac{1}{|\com|}\int_0^{T^0\cos^2\omega} \hat\Gamma_0^{-1}(s)\dd s = R_0 - \frac{1}{|\com|}\int_0^{\infty} z\hat\Gamma_0'(z)\dd z\\[2mm] \label{romb}
&  = R_0+ \frac{1}{|\com|} \int_0^{\infty}\hat\Gamma_0(z)\dd z = R_0 + \frac{u^*}{|\com|} \quad \forall t>0,
\end{align}
while Proposition~\ref{xp3} guarantees only
$$
R^\omega(t) \le \frac{R_0 + u^*}{|\com|} \quad \forall t>0.
$$

We distinguish two different possibilities.

{\bf Case~1}. $P^\flat > 0$. Then necessarily $P(0) > 0$. In this case, a counterpart of \eqref{defq} is meaningful for
\be{defqp}
c > P^\flat \cos\omega = - P^\flat |\cos\omega|.
\ee
For all $z\in (-P^\flat, \infty)$ we have
$$
\frac{Q^*}{(z+P^\sharp)^2} \le \hat\Gamma_0(z) \le \frac{Q^*}{(z+P^\flat)^2},
$$
hence,
$$
\frac{Q^*}{(\hat\Gamma_0^{-1}(\tau)+P^\sharp)^2} \le \tau \le \frac{Q^*}{(\hat\Gamma_0^{-1}(\tau)+P^\flat)^2} \quad \forall \ 0 < \tau < T^\flat\cos^2\omega,
$$
that is,
\be{g0ta}
\sqrt{\frac{Q^*}{\tau}}-P^\sharp \le \hat\Gamma_0^{-1}(\tau) \le \sqrt{\frac{Q^*}{\tau}} - P^\flat \quad \forall \ 0 < \tau < T^\flat\cos^2\omega.
\ee
From \eqref{romdota} and \eqref{g0ta} we thus obtain the estimate
\be{esroa}
R_0 - t\,P^\sharp|\com| + 2\sqrt{Q^* t}  \le R^\omega(t) \le R_0 - t\,P^\flat|\com| + 2\sqrt{Q^* t} \quad \text{for } 0\leq t < T^\flat.
\ee
Recall that at time $t$, the traveling wave $u_\omega(x,t)$ given by \eqref{uom} vanishes for $\scal{x,\bfe_\omega} \ge \rho$ with $\rho$ given by formula \eqref{eom}. Hence, by virtue of \eqref{bom}, \eqref{yan3}, and \eqref{esroa}, we have the implication
\be{bom2a}
\scal{x,\bfe_\omega} > R^\omega_\flat(t)\coloneqq R_0 + t\,P^\flat\com + 2\sqrt{Q^* t} \ \Longrightarrow \ u(x,t) = 0.
\ee
For each fixed time $t$, the wet region $W(t) = \{x \in (-M, M)^N: u(x,t) > 0\}$ is therefore contained in the intersection of half-spaces
\be{bom3a}
W(t) \subset W_-(t)\coloneqq \bigcap_{\omega\in (\pi/2,3\pi/2)} V_\omega(t), \quad V_\omega(t) = \{x \in \real^N: \scal{x,\bfe_\omega} \le R^\omega_\flat(t)\}
\ee
over all $\bfe_\perp$ perpendicular to $\bfe_N$.

\begin{figure}[htb]
\begin{center}
\includegraphics[height=10cm]{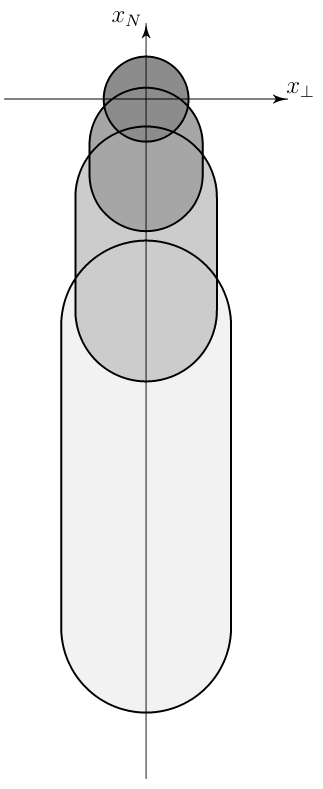}\qquad
\includegraphics[height=10cm]{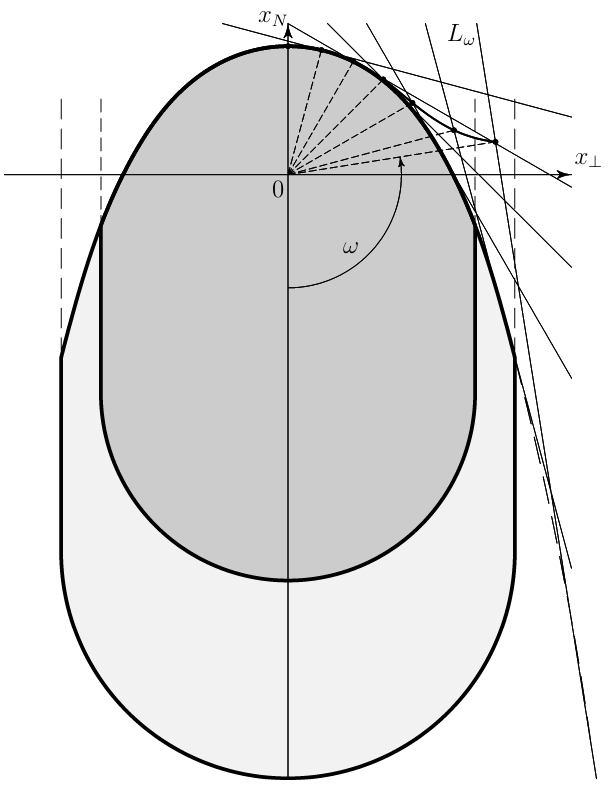}
\caption{Moving wetting fronts intersected with the plane $E(\bfe_\perp, \bfe_N)$ in the case $P(0) > 0$ (left) and $P(0) = 0$ (right).}
\label{f2}
\end{center}
\end{figure}

For $\omega \in (\pi/2,3\pi/2)$ we have $\bfe_{\omega-\pi} = - \bfe_\omega$ and $\cos(\omega-\pi) = -\com$. Hence, from \eqref{bom2}--\eqref{bom3} and \eqref{bom2a}--\eqref{bom3a} we deduce the implication
\be{imp}
u(x,t) > 0 \ \Longrightarrow \ -R_0 + t\,P^\sharp\com - 2\sqrt{Q^* t} < \scal{x,\bfe_\omega} < R_0 + t\,P^\flat\com + 2\sqrt{Q^* t}
\ee
for all $\omega \in (\pi/2,3\pi/2)$. Hence, by the argument leading to \eqref{bound1}, the set $W(t)$ is contained in the convex hull of the union $B_* \cup B^*$ of the two balls
\begin{align*}
B_* &= \{x \in \real^N: |x_\perp|^2 + (x_N + tP^\flat)^2 < (R_0 + 2\sqrt{Q^*t})^2\},\\[2mm]
B^* &= \{x \in \real^N: |x_\perp|^2 + (x_N + tP^\sharp)^2 < (R_0 + 2\sqrt{Q^*t})^2\},
\end{align*}
with $x_N = \scal{x,\bfe_N}$ and $\scal{x_\perp, \bfe_N} = 0$. The time evolution of the wetting front is represented in Figure~\ref{f2} left. Note that, in particular, when $P^\flat = P^\sharp$ the two balls $B_*, B^*$ coincide, and the wet region is contained in the ball of radius $R_0 + 2\sqrt{Q^*t}$ moving downward with speed $P^\flat = P^\sharp$.

{\bf Case~2}. $P^\flat = 0$. Then necessarily $P(0) = 0$ and $P(v) > 0 $ for $v>0$, and \eqref{eotp} is meaningful if and only if $\dot R^\omega(t) \ge 0$. According to the notation in \eqref{tfla}, we distinguish two subcases.

{\bf Case~2a}. $T^\flat\cos^2\omega = \displaystyle{\int_0^{u^*} \frac1{P(v)}\dd v = \infty}$. The function $\hgo$ in \eqref{gamo} is defined for every $z>0$, it is decreasing in $(0,\infty)$, $\lim_{z\to \infty}\hgo(z) = 0$, $\lim_{z\to 0}\hgo(z) = +\infty$, and formula \eqref{romdota} for the moving front holds for every $t \ge 0$. We see in particular that in every direction $\omega\in (\pi/2,3\pi/2)$, the function $R^\omega$ is increasing in the whole time interval $(0,\infty)$. This means that unlike in Case~1, some humidity persists in the region $x_N>0$ for all times.

Let us derive some consequences of \eqref{romdota}. Substituting $s = \hgo(z)$ in \eqref{romdota} and integrating by parts we get
\begin{align*}
R^\omega(t) - R_0 &= \frac{1}{|\com|}\int_0^{t\cos^2\omega} \hgom(s)\dd s
= -\frac{1}{|\com|}\int_{\hgom(t\cos^2\omega)}^\infty z\hgo'(z)\dd z\\
&= t|\com|\, \hgom(t\cos^2\omega) + \frac{1}{|\com|}\int_0^{u^*}\frac{P(v)}{\hgom(t\cos^2\omega) +P(v)}\dd v.
\end{align*}
Passing to the limit as $t\to \infty$ we thus get in addition to \eqref{romb} that
\be{romo}
\lim_{t\to \infty} R^\omega(t) = R_0 + \frac{u^*}{|\com|}.
\ee
In some cases, we can estimate the convergence rate. For example, assuming that $P$ is monotone and
\be{a0}
(P^{-1})'(x) \le \alpha \quad \for x\in (0,u^*)
\ee
for some $\alpha>0$, we get from \eqref{gamo} by substituting $v=P^{-1}(x)$ that
$$
\hgo(z) = \int_0^{P(u^*)} \frac{x (P^{-1})'(x)}{(z+x)^2}\dd x \le \alpha \int_0^{P(u^*)} \frac{1}{z+x}\dd x \le \alpha\log\left(\frac{z+P(u^*)}{z}\right),
$$
and from \eqref{romdo} we conclude that the decay of $\dot R^\omega(t)$ in time is estimated from above as
\be{ese}
\dot R^\omega(t)\le \frac{P(u^*)|\com|}{\expe^{t\cos^2\omega/\alpha}-1},
\ee
hence, it is exponential. Similarly, in the case
\be{a1}
(P^{-1})'(x) \le \alpha x^{-1/m}\quad \for x\in (0,u^*)
\ee
for some $\alpha>0$ and $m>1$, we obtain
\begin{align*}
\hgo(z) &= \int_0^{P(u^*)} \frac{x (P^{-1})'(x)}{(z+x)^2}\dd x \le \alpha\int_0^{P(u^*)} (z+x)^{-1 -(1/m)}\dd x\\
&\le \alpha m\left(z^{-(1/m)}-(P(u^*) + z)^{-(1/m)} \right)
\le \alpha m z^{-1/m},
\end{align*}
hence, $\dot R^\omega(t)$ decays proportionally to $t^{-m}$.

{\bf Case 2b}. $T^\flat\cos^2\omega = \displaystyle{\int_0^{u^*} \frac1{P(v)}\dd v \eqqcolon K < \infty}$. The function $\hgo$ in \eqref{gamo} is defined for every $z>0$, it is decreasing in $(0,\infty)$, $\lim_{z\to \infty}\hgo(z) =0$, and $\lim_{z\to 0}\hgo(z) = K$. As a consequence, the function $\hgom$ is defined only in the interval $(0,K]$, $\lim_{z\to 0}\hgom(z) = +\infty$, $\hgom(K) = 0$, and formula \eqref{romdota} for the moving front holds in $[0,T^\flat)$. By \eqref{yan3}, we can continuously extend the function $R^\omega$ by
\be{romo2}
\dot R^\omega(t) = 0, \quad R^\omega(t) = R_0 + \frac{u^*}{|\com|} \quad \for t \ge T^\flat.
\ee
In other words, the wetting front reaches the upper bound \eqref{romb} asymptotically as $t\to\infty$ in Case~2a, and in finite time in Case~2b.

Overall, since $\dot R^\omega(t) \ge 0$ for all $t > 0$, in Case~2 the wet region $W = \{x \in (-M, M)^N: \exists t\ge 0: u(x,t) > 0\}$ is contained in the intersection of half-spaces
\be{bom3i}
W^\infty = \bigcap_{\omega\in (\pi/2, 3\pi/2)}V_\omega^\infty, \quad V_\omega^\infty = \left\{x \in \real^N: \scal{x,\bfe_\omega} \le  R_0 - \frac{u^*}{\com}\right\}
\ee
over all $\bfe_\perp$ perpendicular to $\bfe_N$. For each choice of $\bfe_\perp$, we estimate from above the intersection $W_{\bfe_\perp}^\infty$ of $W^\infty$ with the plane $E(\bfe_\perp, \bfe_N) \subset \real^N$ generated by $\bfe_\perp$ and $\bfe_N$ with coordinates $x_\perp, x_N$ by the condition that the lines
\be{line}
L_\omega = \left\{x\in E(\bfe_\perp,\bfe_N): \scal{x,\bfe_\omega} = x_\perp \som - x_N \com =  R_0 - \frac{u^*}{\com}\right\}
\ee
are tangent to $W_{\bfe_\perp}^\infty$ for all $\omega \in (\pi/2, 3\pi/2)$. In other words, we determine $W_{\bfe_\perp}^\infty$ as the region bounded by the envelope of the lines $L_\omega$.
Let $x_\perp(\omega)\bfe_\perp + x_N(\omega)\bfe_N$ be the contact point of $L_\omega$ with $W_{\bfe_\perp}^\infty$. Then the tangent vector $x_\perp'(\omega)\bfe_\perp + x_N'(\omega)\bfe_N$ is orthogonal to the normal vector $\som\,\bfe_\perp - \com\,\bfe_N$ of $L_\omega$, that is,
$$
x_\perp'(\omega) \som - x_N'(\omega) \com = 0.
$$
We differentiate the identity in \eqref{line} and get
$$
x_\perp'(\omega) \som - x_N'(\omega) \com + x_\perp(\omega) \com + x_N(\omega) \som = -\frac{u^*\som}{\cos^2\omega}.
$$
Hence, $x_\perp(\omega)$ and $x_N(\omega)$ are determined as solutions to the linear system
\begin{align*}
x_\perp(\omega) \som - x_N(\omega) \com &= R_0- \frac{u^*}{\com},\\
x_\perp(\omega) \com + x_N(\omega) \som &=  -\frac{u^*\som}{\cos^2\omega},
\end{align*}
so that an upper bound for the upward wetting front admits an explicit representation
\be{sol}
x_\perp(\omega) = \left(R_0+ \frac{2u^*}{|\com|}\right)\som, \quad x_N(\omega) = R_0|\com| +u^*\left(1 - \frac{\sin^2\omega}{\cos^2\omega}\right),
\ee
see Figure~\ref{f2} right. Note that a straightforward computation yields that
$$
|x_\perp(\omega)|^2 + |x_N(\omega)|^2 = \left(R_0 + \frac{u^*}{|\com|}\right)^2,
$$
so that formula \eqref{sol} is in agreement with \eqref{romo} and \eqref{romo2}.

We conclude the proof of Theorem~\ref{3t1} by choosing any $M > R_0 + u^* + T P^\sharp + 2 \sqrt{Q^* T}$. Then, according to \eqref{esro} and Proposition~\ref{xp3}, the moisture front stays away from the boundary of $\Omega$ for all times $t\in [0,T]$. We see, in particular, that the solution we have constructed in Theorem~\ref{tfg1} is at the same time a solution to Problem~\eqref{e0u}--\eqref{e2} with homogeneous Neumann or Robin boundary conditions on $\partial\Omega$, and that it can be extended by $0$ to the whole space $\real^N$.


\section{Numerical tests}\label{sec:nume}

\begin{figure}[htb]
	\centering
	\begin{minipage}{0.45\linewidth}
		\centering
	    \includegraphics[width=\linewidth]{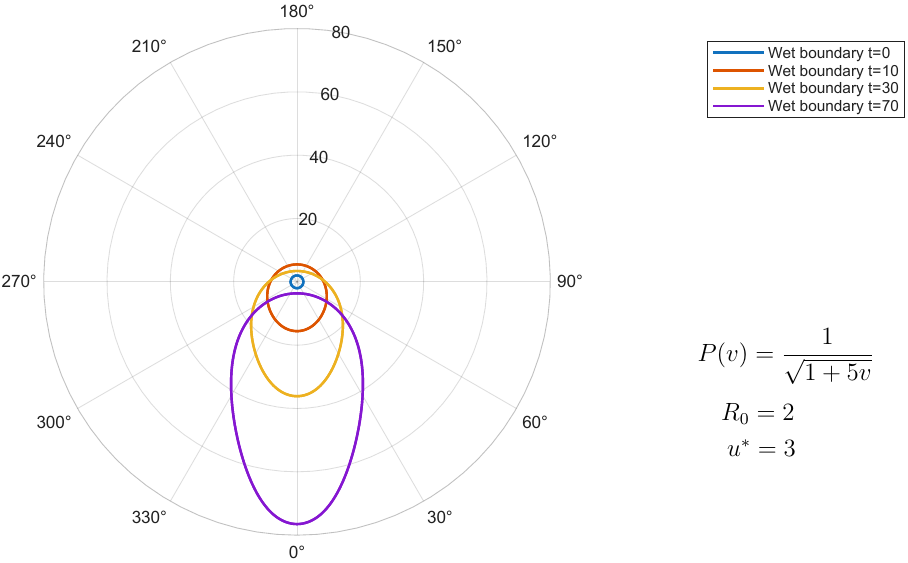}
	    \subcaption{}
	    \label{fig:subfig1}
	\end{minipage}
	\begin{minipage}{0.45\linewidth}
		\centering
	    \includegraphics[width=\linewidth]{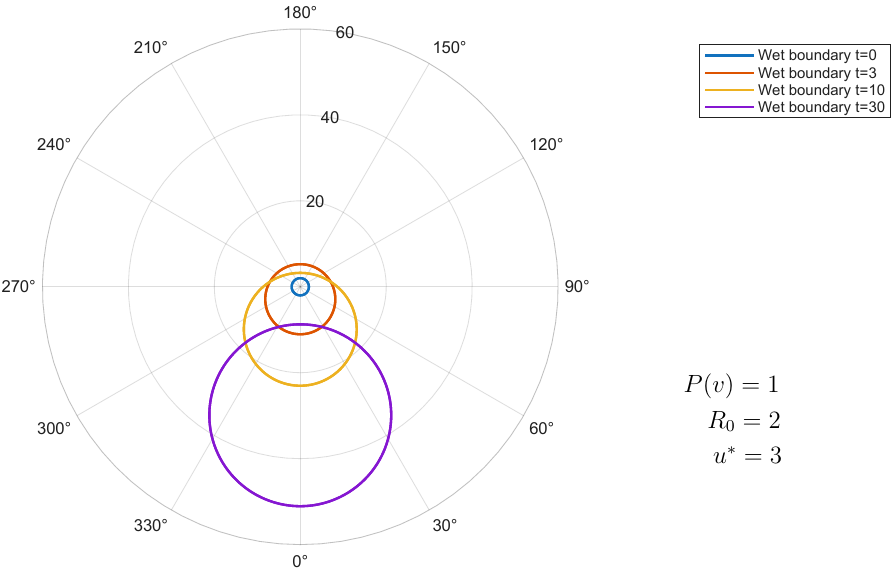}
	    \subcaption{}
	    \label{fig:subfig2}
        \end{minipage} \\	
	\begin{minipage}{0.45\linewidth}
		\centering
	    \includegraphics[width=\linewidth]{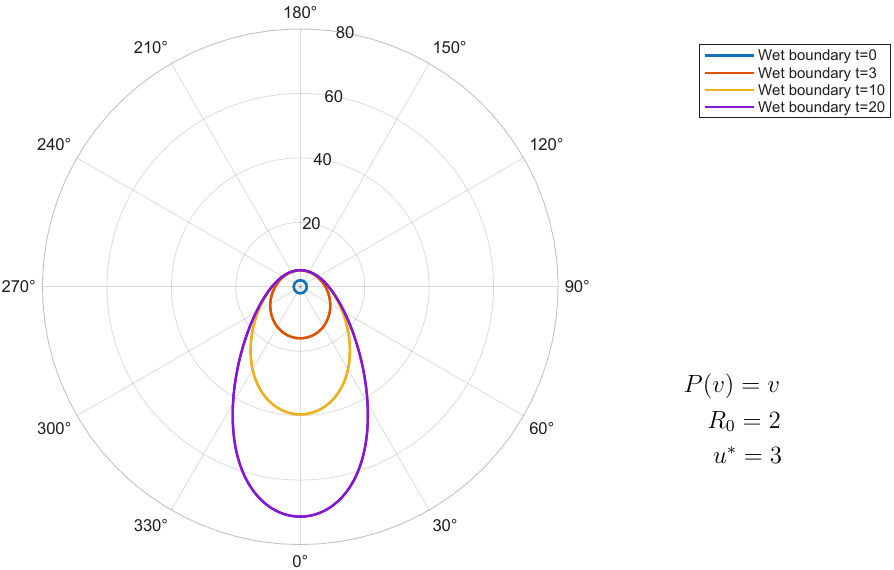}
	    \subcaption{}
	    \label{fig:subfig3}
    \end{minipage}
    \begin{minipage}{0.45\linewidth}
       	\centering
	   	\includegraphics[width=\linewidth]{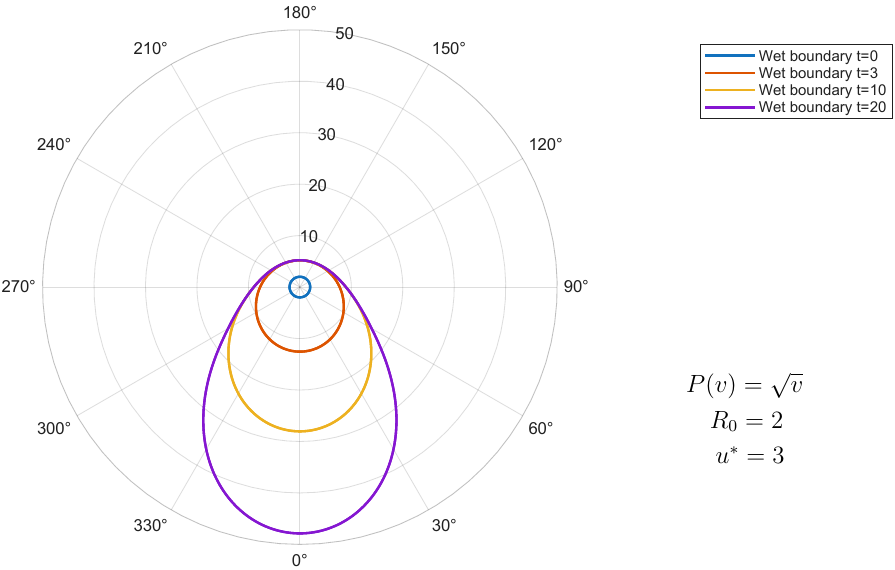}
	   	\subcaption{}
	   	\label{fig:subfig4}
    \end{minipage}
\caption{Moving wetting fronts for different functions $P(v)$. Figures~\ref{fig:subfig1} and \ref{fig:subfig2} illustrate Case~1, while Figures~\ref{fig:subfig3} and \ref{fig:subfig4} illustrate Cases~2a and 2b, respectively.}
\label{num}
\end{figure}

We have shown in Subsection~\ref{shape} that the shape of the wet region is completely determined by the function $P(v)$. In this section, we illustrate the theory developed in Section~\ref{sec:mois} by providing some typical scenarios related to a canonical choice of $\bfe_\perp$ and to different functions $P(v)$, namely, $P(v) = (1 + 5v)^{-1/2}$, $P(v) = 1$, $P(v) = v$, and $P(v) = v^{1/2}$. In Figure~\ref{num} we show the evolution of the moisture front in four consecutive times.

\section{Conclusions}

We have presented a comprehensive study of moisture propagation in homogeneous soils, bridging laboratory tests with rigorous mathematical analysis. The core of our theoretical investigation established that the wetting front exhibits anisotropic behavior: the downward propagation is driven by gravity and scales linearly with time ($t$), while the horizontal spread is dominated by capillary diffusion and scales as $\sqrt{t}$. Furthermore, the analysis correctly captures the phenomenon of upward capillary flow, predicting that the upward front remains bounded, a behavior clearly observed in the experiments where capillary rise arrests after a finite time. We also proved that distinct asymptotic regimes emerge depending on the specific hydraulic properties of the soil. In one regime, the upward front eventually reverses direction, moving downward with a constant asymptotic speed after an initial expansion. In other regimes, the front monotonically approaches its maximum height, reaching it either in finite time or asymptotically as $t \to \infty$. These analytical findings are in agreement with laboratory experiments conducted on fine sand, as the theoretically predicted ``egg-shaped'' geometry confirms the observed wetting bulb. Numerical simulations were employed to illustrate these mathematical results in specific cases. Future work may extend this analytical framework to incorporate hysteresis effects in the soil's hydraulic properties.


\end{document}